\documentclass{amsart}

\usepackage{xypic}

\begin{document}

\newtheorem{proposition}[equation]{Proposition}
\newtheorem{lemma}[equation]{Lemma}
\theoremstyle{definition}
\newtheorem{definition}[equation]{Definition}
\newtheorem{remark}[equation]{Remark}
\newtheorem{example}[equation]{Example}
\theoremstyle{plain}
\numberwithin{equation}{section} 

\title[functional equations of dilogarithms]{describing all multivariable functional equations of dilogarithms}

\author{Rob de Jeu}
\address{Faculteit der B\`etawetenschappen\\Afdeling Wiskunde\\Vrije Universiteit Amsterdam\\De Boelelaan 1111\\1081~HV Amsterdam\\The Netherlands}

\begin{abstract}
We prove quite general statements about functional equations
in any number of variables for the
dilogarithms defined by Bloch-Wigner, Rogers, and Coleman,
showing
that they follow from certain 5-term and 2-term relations in a precise way. Unlike
many other references, we use arbitrary coefficients, do
not ignore any torsion, and get sharp results.
For the Bloch-Wigner dilogarithm, we also consider complex conjugation,
and for Coleman's $ p $-adic dilogarithm we show independence of the branch.
\end{abstract}

\subjclass[2010]{Primary: 11G55; secondary: 19D45}

\keywords{Bloch group, dilogarithm, functional equation, 5-term relation}

\thanks{The author wants to thank Herbert Gangl for useful
comments and discussions, and for the suggestion to write this
paper.}

\maketitle

\def\mybullet{$ \bullet $\hskip5pt}

\def\vph{\vphantom{$b^{b^b}p_{p_{p_p}}$}} 

\def\rightiso{\buildrel{\simeq}\over{\rightarrow}}

\def\mybigcap #1 {{\textstyle \bigcap\limits#1}}

\def\ol{\overline}

\def\C{\mathbb{C}}
\def\P{\mathbb{P}}
\def\Q{\mathbb{Q}}
\def\R{\mathbb{R}}
\def\Z{\mathbb{Z}}

\def\a{\alpha}
\def\ta{\widetilde\a}
\def\b{\beta}
\def\bb{\ta} 
\def\c{\gamma} 

\def\dd{\mathrm{d}}
\def\wirt #1 {\frac{\partial}{\partial #1}} 

\def\id{\mathrm{id}}
\def\im{\mathrm{im}}

\def\ww{\widetilde w} 

\def\tF{\widetilde F}

\def\l{\ol l}
\def\L{\ol{L}}
\def\RL{\ol{RL}}
\def\Li{\mathrm{Li}_2}
\def\Lip{\mathrm{Li}_{p,2}}
\def\Dp{D_p}

\def\sp #1 {\mathrm{sp}_{#1}}
\def\spp #1 {\mathrm{sp}_{\pi_{#1}}}

\def\ord{{\mathrm{ord}}}

\def\Res{{\mathrm{Res}}}

\def\Rft #1 {R_{5,2}(#1)}

\def\Tor{{\rm{Tor}}}
\def\ind{{\mathrm{ind}}}
\def\Kind #1 {K_3(#1)^\ind}
\def\pp{\mathfrak{p}}
\def\pptilde{\widetilde{\pp}}
\def\sustimes #1 #2 {#1 \overset\sigma\otimes #2}
\def\suswedge #1 {(#1 \otimes #1)_\sigma}
\def\tw{\mathrel{\tilde\wedge}}
\def\twt{\tilde\wedge^2}
\def\Btilde{\widetilde B}

\section{Introduction}

Let $ D : \C \setminus \{0,1\} \to \R $ be the Bloch-Wigner dilogarithm,
which is defined in \cite{bl00} by integrating the 1-form $ \log|w| \, \dd \, {\arg}(1-w) - \log|1-w| \, \dd \, {\rm arg}(w) $
along any path from any point in $\R \setminus \{0,1\} $ to $ z $.
The identities
\begin{equation} \label{Dfuns}
\begin{gathered}
 D(z) + D(z^{-1}) = 0,
\qquad
 D(z) + D(1-z) = 0,
\qquad
 D(z) + D(\ol{z}) = 0,
\\
 D(x)-D(y)+D\Bigl(\frac{y}{x}\Bigr)+D\Bigl(\frac{1-x}{1-y}\Bigr)-D\Bigl(\frac{1-x^{-1}}{1-y^{-1}}\Bigr) = 0,
\end{gathered}
\end{equation}
for $ x $, $ y $, and $ z $ in $ \C \setminus \{0,1\} $ with~$ x \ne y $
are easily verified using differentiation.

There is a folklore conjecture that `every functional
equation of~$ D $ can be obtained from 5-term relations',
but this is imprecise, and besides, there are several
closely related but still different 5-term relations. E.g., apart
from the last one in~\eqref{Dfuns}, one also has
$ D(x) + D(y) + D(1-xy)+D(\frac{1-x}{1-xy}) + D(\frac{1-y}{1-xy}) = 0 $
for~$ x, y \ne 0,1 $ and $ xy \ne 1 $.
The statement can be interpreted as follows.
Let~$ F $ be a finitely generated field extension of~$ \C $, which
is the function field of a variety~$ X $. Suppose that for some
$ f_j \ne 0, 1 $ in~$ F $, and integer or complex coefficients $ a_j $,
we have that~$ \sum_j a_j D(f_j) $ is constant 
on the Zariski open part of $ X $ where it is defined.
Subtracting from this the evaluation of this expression at a point of $ X $
where it is defined gives a functional equation, and one expects
it can be expressed as a linear combination of 5-term relations
and perhaps its degeneracies.
The philosophy (or hope) behind such a statement is presumably
based on the known case where~$ F = \C(t) $ for a variable~$ t $ and the $ a_j $ are integers
(see~\cite[Proposition~4(ii)]{Zag07} and cf.~\cite[\S4]{Woj91}),
and is strengthened by another result in this direction (see \cite[\S5]{souderes18}).

We describe all such functional equations with coefficients~$ a_j $
in a subgroup of $ \C $ if $ F = k(t_1,\dots,t_n) $
where~$ k $ is a subfield of $ \C $ and the $ t_j $ are variables
(see Proposition~\ref{funprop} and the discussion preceding it),
and similarly for Rogers's dilogarithm (see Proposition~\ref{Rfunprop})
and Coleman's $ p $-adic dilogarithm (see Proposition~\ref{Coleprop}).
For the Bloch-Wigner dilogarithm we also consider functional
equations involving conjugation (see Propositions~\ref{conjpropR} and~\ref{conjpropC}).
All such functional equations can be obtained in a precise way from the first and
last functional equations given in~\eqref{Dfuns} or their equivalents
for the Rogers and Coleman dilogarithms.
The proofs, although they share a common pattern, are a bit technical
and we postpone them until Section~\ref{proofs}, opting to give
full statements of our results
in Section~\ref{results} first.
We also postpone some remarks on sharpness of our results
(see Remarks~\ref{sharpDrem} and~\ref{Dbarsharprem}.)

Our approach to proving Proposition~\ref{funprop} for one variable is very
different from that used to prove~\cite[Proposition~4(ii)]{Zag07}.
The latter is based on rewriting the arguments in a functional equation
using those listed in~\eqref{Dfuns}, which is a modification
of an argument in~\cite[\S4]{Woj91} and requires the base field
to be algebraically closed. Instead,
working directly over the base field~$ k $,
we translate any functional equation into a boundary condition in a suitable wedge product (see below)
by means of differentiation
(cf.~ the proofs of~\cite[Proposition~4.9, p.79]{gonXpam} and~\cite[Proposition~4(ii)]{Zag07}).
For two or more variables we deduce the boundary condition
from the case of one variable by using Lemma~\ref{buddinglemma}.
Finally, we specialise the variables, and use a modification of a
(seemingly little known)
result by Suslin (see Propositions~\ref{susprop} and~\ref{barsusprop}) repeatedly to finish the proof.

This method also applies to Propositions~\ref{conjpropR}, \ref{Rfunprop} and~\ref{Coleprop},
with the details in the case of one variable depending on the dilogarithm involved.
By contrast, Proposition~\ref{conjpropC} is
deduced from Proposition~\ref{conjpropR}.

We finish this introduction by mentioning that
for our fields $ F = k(t_1,\dots,t_n) $
the result of Suslin that we use
is equivalent to the natural map
$ \Kind k \to \Kind F $ being an isomorphism (see the proof
of Proposition~\ref{susprop}).
If $ k $ is the algebraic closure of itself in the function field~$ F $ of
a variety over~$ k $, then such an isomorphism would
presumably make it possible to describe all functional equations of the dilogarithms with arguments in $ F $.
The isomorphism would be a consequence of the `rigidity conjecture', which states
that for any field~$ F $ the natural map $ \Kind F_0 \to \Kind F $ is an isomorphism
if $ F_0 $ is the algebraic closure in $ F $ of its prime field
(see \cite[5.4]{sus:aka} or \cite[VI.5.3.1]{WeiKbook}).
However, the rigidity conjecture, or even this consequence, currently seems well out of reach.

\section{The results} \label{results}

For any field $ F $, we let $ F^\flat = F \setminus \{0,1\} $,
denote by $ \Z[F^\flat] $ the free Abelian group on $ F^\flat $,
and for any torsion-free additive  group $ A $ set~$ A[F^\flat] = A \otimes_\Z \Z[F^\flat] $.
We also define
\begin{equation*}
R_5(F) = \Bigl\langle [x]-[y]+\Bigl[\frac{y}{x}\Bigr]+\Bigl[\frac{1-x}{1-y}\Bigr] - \Bigl[\frac{1-x^{-1}}{1-y^{-1}}\Bigr] \text{ with } x, y \text{ in } F^\flat, \ x \ne y \Bigr\rangle
\,,
\end{equation*}
the subgroup of $ \Z[F^\flat] $ generated by the 5-term relations
corresponding to one of the functional equations in~\eqref{Dfuns},
and write $ A R_5(F) $ for the subgroup $ A \otimes_\Z R_5(F) $
of~$ A[F^\flat] $.
(Every finitely generated subgroup of~$ A $ is isomorphic to $ \Z^l $
for $ l \ge 0 $, and~$ A $ is the union of such subgroups, so tensoring with $ A $ over $ \Z $ is exact.)
Similarly, we introduce the subgroup
\begin{equation*}
R_2(F) = \langle [x]+[x^{-1}] \text{ with } x \text{ in } F^\flat \rangle
\end{equation*}
of $ \Z[F^\flat] $ and the subgroup $ A R_2(F) $
of $ A[F^\flat] $, which corresponds to another functional equation
in~\eqref{Dfuns}. We note that $ 2 R_2(F) \subseteq R_5(F) $ if $ |F| \ge 4 $
by~\cite[Lemma~1.2]{susXkof} or~\cite[VI.5.4]{WeiKbook}
but we shall obtain better results by not using this.
We shall write $ \Rft F $ for $ R_5(F) + R_2(F) $ and $ A \Rft F $ for $ A R_5(F) + A R_2(F) $.

If $ t $ is a variable and $ b $ is in $ F \cup \{ \infty \} $, then
one can specialise $ t $ to $ b $, so
any~$ f(t) $ in~$ F(t)^\flat $ evaluates to an element of $ F \cup \{\infty\} $,
giving a map $ \Z[F(t)^\flat] \to \Z[F \cup \{\infty\}] $.
It is customary to map the latter to $ \Z[F^\flat] $ by mapping
$ [0] $, $ [1] $ and $[\infty] $ to~$ 0 $, but we shall use a
different convention so discuss this in some detail.

The result of specialising $ t $ to $ b $ in $ [x]-[y]+\bigl[\frac{y}{x}\bigr]+\bigl[\frac{1-x}{1-y}\bigr] - \bigl[\frac{1-x^{-1}}{1-y^{-1}}\bigr] $
with $ x \ne y $ in $ F(t)^\flat $ is as follows, where~$ x(b) $
is given horizontally and $ y(b) $ vertically, and $ c $ is some
element of $ F^\flat $.
The entries are easily checked by some mild rewriting.
E.g., if~$ x(b) = y(b) = \infty $ then use $ \frac{1-x}{1-y} = \frac xy\frac{1-x^{-1}}{1-y^{-1}} $.
\begin{center}
\smaller\smaller
\begin{tabular}{|c|c|c|c|c|}
\cline{1-5}
\vph & $ 0 $ & $ 1 $ & $  \infty $ & other
\\
\cline{1-5}
\vph
$ 0 $ & $ [1] $ & $ [1] $ & $ 2 [\infty] - [0] $ & $ [c] + [1-c] - [0] $
\\
\cline{1-5}
\vph
$ 1 $ & $  [0] - [1] + [\infty] $ & $ [1] $ & $ [\infty] - [1] + [0] $ & $ [c] - [1] + [c^{-1}]  $
\\
\cline{1-5}
 \vph
 \lower-7pt\hbox{$ \infty $} & \lower-7pt\hbox{$ 2 [0] - [\infty] $} & \lower-7pt\hbox{$ [1] $} & \vbox{\hbox{\vph $ [1] $ or $ [0] + [\infty] - [1] $}\hbox{or $ [c] + [c^{-1}] - [1] $}} & \lower-7pt\hbox{$ [c] + [0] - [1-c^{-1}] $}
\\
\cline{1-5}
\vph
other & $ [0] - [c] + [(1-c)^{-1}] $ & $ [1] $ & $  2 [\infty] + [0] - [c] - [(1-c^{-1})^{-1}] $ & $ [1] $ or in $ R_5(F) $
\\
\cline{1-5}
\end{tabular}
\end{center}

\smallskip
\noindent
Similarly, $ [x] + [x^{-1}] $ can specialise to $ [c] + [c^{-1}] $
for $ c $ in $ F^\flat $, $ 2[1] $ or $ [0] + [\infty] $.
Therefore $ \Rft F(t) $ is mapped to the subgroup
\begin{equation*}
 \Rft F + \langle [1] , [0] + [\infty] , 3 [0] \rangle + \langle [c] + [1-c] - [0]
\text{ with } c \text{ in } F^\flat \rangle
\end{equation*}
of $ \Z[F \cup \{\infty\} ] = \Z[F^\flat] \oplus \Z[\{0, 1, \infty\}] $.

Now assume that $ |F| \ge 4 $.
If we fix $ c $ in $ F^\flat $, then $ C_{F,c} = [c] + [1-c] $ in $ \Z[F^\flat] $
is such that its class modulo~$ R_5(F) $ is independent of $ c $
(see~\cite[Lemma~1.3]{susXkof} or~\cite[VI.5.4]{WeiKbook}),
so we can rewrite this subgroup as
\begin{equation}  \label{specsubgroup}
 \Rft F + \langle [1] , [0] + [\infty] , 3 [0] \rangle + \langle C_{F,c} - [0] \rangle
\,.
\end{equation}
We then define $ \sp b,c $, the specialisation of $ t $ to $ b $
for this fixed $ c $, as the composition
\begin{equation*}
\sp b,c : \Z[F(t)^\flat] \to \Z[F \cup \{\infty\}] \to \Z[F^\flat] 
\end{equation*}
with the first map induced by evaluation of functions at $ b $,
and the second given by
\begin{equation*}
 \a  \mapsto \a - c_1(\a) [1] + c_0(\a) (C_{F,c} - [0]) - c_\infty(\a) (C_{F,c} + [\infty]) 
\end{equation*}
with~$ c_x(\a) $ the coefficient of $ [x] $ in $ \a $.
Clearly~$ \sp b,c $ maps~$ R_2(F(t)) $ to $ R_2(F) $. It maps~$ \Rft F(t) $ to~$ \Rft F $ because
$ 3 C_{F,c} $ is in $ \Rft F $
(see~\cite[Lemma~1.4]{susXkof} or~\cite[VI.5.4]{WeiKbook}).
Note that we chose the second map such that it induces an isomorphism
\begin{equation} \label{introiso}
\frac{\Z[F^\flat] \oplus \Z[\{0, 1, \infty\}]}{\Rft F + \langle [1] , [0] + [\infty] , 3 [0] \rangle + \langle C_{F,c} - [0] \rangle}
\rightiso
\frac{\Z[F^\flat]}{\Rft F }
\end{equation}
that is independent of the choice of $ c $, and is, in
fact, the inverse of the natural map in the other direction.

The corresponding statements with coefficients~$ A $ as before also hold.
For simplicity, we use the notation $ \sp b,c $ for what technically
should be $ \id_A \otimes \sp b,c $.

With more variables it is more complicated.
For example, for~$ F(t_1,t_2) $ we can first specialise~$ t_1 $ to
$ b $ in $ F(t_2) \cup \{\infty\} $ (choosing a $ c $ in $ F(t_1)^\flat $),
and then $ t_2 $ in the result to $ b' $ in~$ F \cup \{\infty\} $
(using a $ c' $ in $ F^\flat $).
In general, for~$ F(t_1,\dots,t_n) $ we can specialise 
 $ t_i $ to an element $ b $ of~$ F(t_1,\dots,t_{i-1}, t_{i+1},\dots, t_n) \cup \{\infty\} $
under the choice of a $ c $ in $  F(t_1,\dots,t_{i-1}, t_{i+1},\dots,t_n)^\flat $,
and then apply such a procedure again to the result.
In total we can perform $ n $ such (one step) specialisations
as we run out of variables, but we can also choose to perform fewer
and keep some of the variables in the result.
If we apply this to an element $ \a $ of~$ A[F(t_1,\dots,t_n)^\flat] $
then we call the result an \emph{iterated specialisation} of
$ \a $.

Even if we specialise all $ t_i $ to elements of
$ F \cup \{\infty\} $ the result may depend on the order,
as, e.g., specialising~$ t_1 $ and $ t_2 $ to 0 in 
$ [(t_1 + c t_2)/(t_1 + t_2)] $ for $ c $ in $ F^\flat $ results in
either $ [c] $ or $ 0 $.
But if $ P = (b_1,\dots,b_n) $ in $ F^n $ is such that for each
$ [f] $ occurring in a given $ \a $ 
in $ A[F(t_1,\dots,t_n)^\flat] $ the function $ f $ is defined at~$ P $
with~$ f(P) \ne 0 , 1 $, then for $ \a $ the order does
not matter as the result is the same as simply replacing all
or some of the $ t_i $ with~$ b_i $, and we do not even have
to deal with~$ [0] $, $ [1] $ or $ [\infty] $.
In such a situation we can safely write~$ \a(b_1,\dots,b_n) $ for the result if we specialise each
$ t_i $ to $ b_i $, and use similar notation if we specialise only some of the $ t_i $
to $ b_i $.

\begin{example} \label{introexample}
Since we work at the level of generators and not modulo~$ R_5 $,
the results of iterated specialisation can be surprising
because of how we deal with~$ [0] $ and $ [\infty] $ in
the naive specialisation and the choice of the auxiliary $ c $
involved.
For~$ \a = [t^2] - 2[t] - 2 [-t] $ in~$ \Z[F(t)^\flat] $ we find $ \sp 1,c (\a) = -2 [-1] $
but~$ \sp \infty,c (\a) = 3 [c] + 3 [1-c] $.
For $ \a = [t_1+t_2+t_3] $ in $ \Z[F(t_1,t_2,t_3)^\flat] $, and
specialising~$ t_2 $ to $ -t_1 - t_3 $ with the choice
$ c = t_1 + t_3^2 $ in $ F(t_1,t_3)^\flat $, we find
$ [t_1 + t_3^2] + [1 - t_1 - t_3^2] $ in~$ \Z[F(t_1,t_3)^\flat] $.
Doing a further specialisation, with $ t_1 $ specialising to $ - t_3^2 $
and choosing~$ c' = t_3 $, we obtain~$ [t_3] + [1-t_3] $ in~$ \Z[F(t_3)^\flat] $ as
an iterated specialisation of $ \a $.
\end{example}

Before we state our results in detail, we emphasise an important consequence
for functional equations.
If $ A $ is a subgroup of $ \C $, $ k $ a subfield of $ \C $,
and~$ \a = \sum_j a_j [f_j] $ is in $ A[k(t_1,\dots,t_n)^\flat] $,
then we write~$ D(\a) $ for~$ \sum_j a_j D(f_j) $ and view this as a function
on a suitable Zariski open part of~$ \C^n $.
(We shall tacitly use similar notation for the Rogers and $ p $-adic dilogarithms later.)
If $ \a $ is in $ A[k^\flat] + A \Rft k(t_1,\dots,t_n) $ then~$ D(\a) $ is constant. Conversely, if $ D(\a) $ is constant then
by Proposition~\ref{funprop} there is an~$ \a' $ in~$ A[k^\flat] $
(e.g.,~$ \a(P) $ for a suitable $ P $ in $ k^n $) with $ \a - \a' $
in~$ A \Rft k(t_1,\dots,t_n) $.
So the proposition implies that
\begin{equation*}
\begin{gathered}
\text{\emph{for $ \a $ in $ A[k(t_1,\dots,t_n)^\flat] $, the function $ D(\a) $ is constant}}\\\text{\emph{if and only if $ \a $ is in $ A[k^\flat]  + A \Rft k(t_1,\dots,t_n) $},}
\end{gathered}
\end{equation*}
which describes precisely all functional
equations of $ D $ as discussed before. But the
proposition is more precise because, if $ D(\a) $ is constant,
then it specifies one can use $ \a' = \a(P) $ for a suitable~$ P $ but it also treats any iterated specialisation $ \a' $
of~$ \a $.

Similar descriptions of functional equations for $ D $ in other
situations follow from Propositions~\ref{conjpropR} and~\ref{conjpropC}
using also $ D(z) = - D(\ol z) $ (see Remark~\ref{Dbarsharprem}),
and for the dilogarithms of Rogers and Coleman they follow as
above from Propositions~\ref{Rfunprop} and~\ref{Coleprop}.

The propositions also gives information on $ \a $ in terms its
image in a wedge product, and it is using this that we can
deal with more than one variable in a straightforward
way (see Lemma~\ref{buddinglemma}).
For any field $ F $, we let
\begin{equation*}
 \twt F^* = F^* \otimes_\Z F^* / \langle (-x) \otimes x \text{ with } x \text{ in } F^* \rangle 
\end{equation*}
and write $ x \tw y $ for the class of $ x \otimes y $.
We define $ \partial_F : \Z[F^\flat] \to \twt F^* $ by mapping~$ [x] $ to $ x \tw (1-x) $. It is easy to check 
that $ \Rft F $ is in the kernel of $ \partial_F $.
The inclusion~$ F^* \to F(t)^* $ is split by mapping every
monic irreducible in~$ F[t] $ to~$ 1 $, so we
may view~$ A \otimes_\Z \twt F^* $ as a subgroup of $ A \otimes_\Z \twt F(t_1,\dots,t_n)^* $,
and we shall therefore often write $ \partial $ for either~$ \partial_F $ or~$ \partial_{F(t_1,\dots,t_n)} $.
Moreover, it follows from~\eqref{spCD}
that if $ \a $ in~$ A[F(t_1,\dots,t_n)^\flat] $ is such that $ \partial (\a) $ is
in~$ A \otimes_\Z \twt F^* $, then~$ \partial(\a) = \partial (\a') $
for any iterated specialisation $ \a' $ of $ \a $.

\begin{proposition} \label{funprop}
Let $ A $ be a subgroup of $ \C $, $ k $ a subfield of $ \C $, and let $ t_1,\dots,t_n $
be variables. If $ \a $ in~$ A[k(t_1,\dots,t_n)^\flat] $ is
such that $ D(\a) $ is constant on the Zariski open part of $ \C^n $ where it is defined, then~$ \partial (\a) $ is in~$ A \otimes_\Z \twt k^* $, and~$ \a - \a' $ is
in~$ A R_{5,2}(k(t_1,\dots,t_n)) $
for any iterated specialisation $ \a' $ of~$ \a $.
In particular, if~$ P $ in~$ k^n $ is such
that all functions in $ \a $ are defined at $ P $ with value not equal to 0 or~1,
this holds for~$ \a' = \a(P) $.
\end{proposition}

\begin{remark}
Perhaps surprisingly,
the functional equation $ D(z) + D(1-z) = 0 $
plays no role in this proposition or the ones below.
For simplicity of notation, we discuss the case $ A = \Z $.
Our definition of $\sp b,c $ is such that it removes any
term~$ [c] + [1-c] $ in~\eqref{specsubgroup} that arises by
using a naive (one step) specialisation of any~$ \a $ in $ R_{5,2} $.
Also, if we were to start with~$ \a = [f] + [1-f] $ for $ f $ in $ k(t_1,\dots,t_n)^\flat $
then a (one step) specialisation $ \ta = \sp b,c (\a) $ of $ \a $ is of the form~$ [g] + [1-g] $ or~$ -2 ([c]+ [1-c]) $
as the naive specialisation~$ 2 [\infty] $ is changed into $ -2 ([c]+ [1-c]) $
and~$ [0] + [1] $ into $ [c] + [1-c] $, so $ \a - \ta $ is of the form~$[f] + [1-f] - [g] - [1-g] $ or~$ [f] + [1-f] + 2 [c] +2 [1-c] $. But as mentioned before,
the element~$ [x] + [1-x] $ modulo $ R_5 $ is independent of~$ x $, and $ 3 [x] + 3 [1-x] $
is in~$ R_{5,2} $.
Hence $ \a - \ta $ is in~$ R_{5,2} $.
This also applies to the differences at the remaining steps of the iterated
specialisation, and hence to $ \a - \a' $ by taking a telescoping
sum.

That this functional equation plays no role is no coincidence, because it
does not hold for the version of the Rogers dilogarithm that
vanishes on $ R_{5,2} $ (see Section~\ref{proofs}), and we
adjusted our specialisation to this.
\end{remark}

We now consider statements that also involve complex conjugation,
corresponding to the functional equation~$ D(z) + D(\ol z ) = 0 $.
Because $ \C(z , \ol z) = \C(t_1,t_2) $
where we only plug in real $ t_1 $ and $ t_2 $
(as real and imaginary parts of $ z $), we first consider~$ k(t_1,\dots,t_n) $ with $ k $
a subfield of $ \C $, and $ t_1,\dots,t_n $ taking values in
$ \R $.
Assuming~$ k $ is stable under complex conjugation,
for~$ f $ in~$ k(t_1,\dots,t_n) $ we let~$ \ol f $ be the result
obtained by letting complex conjugation act on the coefficients
of~$ f $, and we use the same notation for elements of $ A[k(t_1,\dots,t_n)^\flat] $
or~$ A \otimes_\Z \twt k(t_1,\dots,t_n)^* $
if~$ A $ is a subgroup of~$ \C $. (Note that we do not let complex conjugation 
act on $ A $.)

\begin{proposition} \label{conjpropR}
Let $ A $ be a subgroup of $ \C $, $ k $ a subfield of $ \C $
that is stable under complex conjugation,
and let $ t_1, \dots, t_n $ be variables.
If~$ \a $ in~$ A[k(t_1,\dots,t_n)^\flat] $ 
is such that~$ D(\a) $ is constant on the Zariski open part of~$ \R^n $ where it is defined, then~$ \partial(\a) $
is the sum of an element in
$ A \otimes_\Z \twt k(t_1,\dots,t_n)^* $ that is invariant under
complex conjugation, and an element of $ A \otimes_\Z \twt k^* $.
The difference between $ \a - \ol \a $ and
any of its iterated specialisations
$($e.g., the evaluation at a point in~$ k^n $ where
all functions in $ \a $ and $ \ol \a $ are defined with value not equal to 0 or~1$)$
is in $ \Rft k(t_1,\dots,t_n) $.
\end{proposition}

\begin{remark} \label{conjremR}
The iterated specialisation of $ \a - \ol \a $ is not necessarily anti-invariant under complex
conjugation. But it is if at every step one uses $ \sp b,c $ with $ b $ and $ c $ invariant
under complex conjugation as it commutes with complex conjugation.
\end{remark}

In the next proposition we consider functions in complex variables and their conjugates,
for which we write~$ z_1, \ol z_1,\dots, z_n, \ol z_n $.
The functions given by polynomial expressions have the structure of a polynomial ring in
$ 2n $ variables (see Section~\ref{proofs}), with a matching
field of fractions, but as we
can specialise only $ n $ variables independently we
limit the type of specialisation.
In the statement below complex conjugation acts on the variables
and on the subfield $ k $ of $ \C $, but
again not on the subgroup $ A $ of $ \C $.

\begin{proposition} \label{conjpropC}
Let $ A $ be a subgroup of $ \C $, $ k $ a subfield of $ \C $
that is stable under complex conjugation,
and let~$ z_1, \dots, z_n $ be variables.
If~$ \a $ in $ A[k(z_1,\ol z_1,\dots,z_n, \ol z_n)^\flat] $
is such that $ D(\a) $ is constant on the Zariski open part of~$ \C^n $ where it is defined,
then~$ \partial(\a - \ol \a) $ is
in $ A \otimes_\Z \twt k^* $.
For any point $ P $ where all the functions involved in~$ \a $
are defined with value not equal to 0 or 1,
we have that $ \a - \ol \a - \a(P) + \ol{\a(P)} $ is
in~$ \Rft k(z_1,\ol z_1,\dots,z_n,\ol z_n $.
\end{proposition}

We now introduce the Rogers dilogarithm~$ L $ on $ \R^\flat $,
for which we follow \cite[p.23]{Zag07}.
With $ \Li(x) = \sum_{n=1}^\infty x^n/n^2 $ for~$ |x| \le 1 $,
we let
$ L(x) = \Li(x) + \frac12 \log(x) \log(1-x) $ for $ 0 < x < 1 $,
$ L(x) = \frac{\pi^2}3 - L(x^{-1}) $ for~$ x >1 $ and
$ L(x) = - L(1-(1-x)^{-1}) $ if~$ x < 0 $.
Then
$ \dd L(x) = \frac12 \log|x| \, \dd \log|1-x| - \frac12 \log|1-x| \, \dd\log|x| $ for
all~$ x $ in~$ \R^\flat $.
It can be extended to a continuous function on $ \R $ by $ L(0) = 0 $
and $ L(1) = \frac{\pi^2}6 $. 
We can view the values in $ \R/\frac{\pi^2}2 \Z $, giving a function~$ \L $
with better functional equations, and which extends to a continuous
function on $ \P_\R^1 $ by $ \L(\infty) = - \l $
with $ \l $ the class of $ \frac{\pi^2}6 $.
The map it induces does not vanish 
on the subgroup in~\eqref{specsubgroup},
but using the function $ \RL(x) = \L(x) - \l $ instead
(cf.~\cite[Definition~2.1]{licht1989},
\cite[p.188]{parry-sah-1983},
\cite[p.646]{frenkel-szenes-1993})
does give such a map (see Section~\ref{proofs}).
Note that in the next proposition the values are in $ \R A / \frac{\pi^2}2 A \subseteq \C / \frac{\pi^2}2 A $.

\begin{proposition} \label{Rfunprop}
Let $ A $ be a subgroup of $ \C $, $ k $ a subfield of $ \R $,
and let $ t_1, \dots, t_n $ be variables.
If~$ \a $ in~$ A[k(t_1,\dots,t_n)^\flat] $ 
is such that~$ \RL(\a) $ is constant on the Zariski open part of~$ \R^n $ where it is defined, then~$ \partial(\a) $
is in $ A \otimes_\Z \twt k^* $, and~$ \a - \a' $ is
in~$ A \Rft k(t_1,\dots,t_n) $
for any iterated specialisation $ \a' $ of~$ \a $.
In particular, if~$ P $ in~$ k^n $ is such
that all functions in $ \a $ are defined at $ P $ with value not equal to 0 or~1,
this holds for $ \a' = \a(P) $.
\end{proposition}

We now move on to Coleman's $ p $-adic dilogarithm.
For this, we fix a branch of the~$ p $-adic logarithm, i.e., a homomorphism
$ \log_p : \C_p^* \to \C_p $ that around 1 is given by the
usual power series. Such a branch is determined by choosing~$ \log_p(p) $ in
$ \C_p $.
Letting~$ \Lip(z) = \sum_{n=1}^\infty z^n/n^2 $ for $ z $ in $ \C_p $
for $ |z|_p < 1 $, Coleman \cite{Col82} uses the action of Frobenius
to extend~$ \Dp(z) = \Lip(z) + \frac12 \log_p(z) \log_p(1-z) $
to a function on $ \C_p^\flat $ that around every point
can be written as a convergent power series, but its
values may depend on the chosen branch of $ \log_p $.
(The precise dependency on the branch is given in the proof of Proposition~\ref{Coleprop}
in Section~\ref{results}.)
One has~$ \dd \Dp(z) = \frac12 \log_p(z) \, \dd \log_p(1-z) - \frac12 \log_p(1-z) \, \dd \log_p(z) $
on $ \C_p^\flat $,
and $ \Dp(z) $ satisfies the functional equations
corresponding to $ R_5 $ and $ R_2 $ (see Section~\ref{proofs}).

\begin{proposition} \label{Coleprop}
Let $ A $ be a subgroup of $ \C_p $, $ k $ a subfield of\,~$ \C_p $,
and let $ t_1, \dots, t_n $ be variables.
If~$ \a $ in~$ A[k(t_1,\dots,t_n)^\flat] $ 
is such that~$ \Dp(\a) $ is constant on the Zariski open part of~$ \C_p^n $ where it is defined, then~$ \partial(\a) $
is in $ A \otimes_\Z \twt k^* $, and~$ \a - \a' $
is in~$ A \Rft k(t_1,\dots,t_n) $ for any
iterated specialisation $ \a' $ of~$ \a $.
In particular, if~$ P $ in~$ k^n $ is such
that all functions in $ \a $ are defined at $ P $ with value not equal to 0 or~1,
this holds for $ \a' = \a(P) $.
Moreover, if $ \Dp^\circ $ is the result of using the branch
$ \log_p^\circ $ of the logarithm, for any~$ \a $ in~$ A[k(t_1,\dots,t_n)^\flat] $,
the difference $ \Dp^\circ(\a) - \Dp(\a) $
can be computed from $ \partial(\a) $ by mapping $ a \otimes ( f \tw g) $ to
\begin{equation*}
 a \frac{\log_p^\circ(p) - \log_p(p)}2 \left( v_p(f) \log_p(g) - v_p(g) \log_p(f) \right)
 \,,
\end{equation*}
for $ v_p $ the valuation on $ \C_p^* $ with $ v_p(p) = 1 $,
where $ v_p(f) \log_p(g) - v_p(g) \log_p(f) $ is independent of the branch of the logarithm used in
it.
In particular, if $ \Dp(\a) $ is constant for one branch, then
it is constant for every branch, and this gives an explicit 
method of computing the difference between the constants.
\end{proposition}

\section{The proofs} \label{proofs}

We begin with recalling various definitions and results from \cite{susXkof}
for infinite fields, which were extended to fields with at least~4 elements
in~\cite[VI.5]{WeiKbook}.

For any field $ F $, we define the pre-Bloch group~$ \pp(F) = \Z[F^\flat] / R_5(F) $
and
\begin{alignat*}{1}
\delta_F : \Z[F^\flat] & \to \suswedge F^*
\\
[x] & \mapsto \sustimes x (1-x)
\,,
\end{alignat*}
where
\begin{equation*}
\suswedge F^* =  \frac{F^*\otimes_\Z F^*}{\langle x \otimes y + y \otimes x \text{ with } x, y \text{ in } F^* \rangle}
\,,
\end{equation*}
and we write $ \sustimes y z $ for the class in the quotient of an element $ y \otimes z $.
As~$ \delta_F (R_5(F)) = 0 $, we can define
the Bloch group $ B(F) = \ker(\delta_F)/R_5(F) $ of $ F $ inside~$ \pp(F) $.

Suslin was interested in its relation to $ \Kind F = K_3(F) / K_3^M(F) $
with $ K_3^M(F) $ the third Milnor $ K $-group of $ F $, which
injects into $ K_3(F) $ (see~\cite[VI.4.3.2]{WeiKbook}).
He proved \cite[Theorem~5.2]{susXkof} that there is a short exact sequence,
natural in~$ F $,
\begin{equation} \label{sussos}
 0 \to \Tor(F^*, F^*)^\sim \to \Kind F \to B(F) \to 0
\end{equation}
if $ F $ is infinite, with
$ \Tor(F^*, F^*)^\sim $ the unique non-trivial
extension of $ \Tor(F^*, F^*) $ by~$ \Z/2\Z $ if the characteristic
of $ F $ is not~2 and $ \Tor(F^*, F^*) $  if it is.
Weibel showed~\eqref{sussos} also exists if $ |F| \ge 4 $ (see~\cite[VI.5.2]{WeiKbook}).

Our main interest is the following consequence.

\begin{proposition} \label{susprop}
If $ |F| \ge 4 $, and $ t $  is a variable,
then the map $ B(F) \to B(F(t)) $ is an isomorphism.
\end{proposition}

\begin{proof}
This is stated (without proof) as Corollary~5.6 in~\cite{susXkof} if $ F $ is infinite.
It can be obtained by mapping the sequence~\eqref{sussos} for
$ F $ with $ |F| \ge 4 $ to that of $ F(t) $, and using 
$ \Tor(F^*,F^*) = \Tor(F(t)^*, F(t)^*) $ as well as~$ \Kind F \rightiso \Kind F(t) $
for any field $ F $.
The latter isomorphism follows from mapping
the short exact sequence
$ 0 \to K_3^M(F) \to  K_3^M(F(t)) \to \coprod_P K_2^M(F[t]/P) \to 0  $
(see \cite[V6.7.1]{WeiKbook})
to the short exact sequence
$ 0 \to K_3(F) \to  K_3(F(t)) \to \coprod_P K_2(F[t]/P) \to 0  $
(see \cite[III.7.4]{WeiKbook})
because $ K_2^M $ and $ K_2 $ of a field coincide.
In both sequences $ P $ runs through the maximal ideals of~$ F[t] $.
\end{proof}

It follows from Proposition~\ref{susprop} that for $ |F| \ge 4 $
and $ \a(t) $ in the kernel of $ \delta_{F(t)}$,
there exists $ \bb $ in the kernel of $ \delta_F $
with~$ \a(t) - \bb $ in $ R_5(F(t)) $.
If $ F $ is infinite then for a given $ \a $ there exist $ b $ in $ F^\flat $
such that evaluating everything at $ t = b $ shows~$ \a(b) - \bb $
is in $ R_5(F) $, so the class of $ \a(b) $ in $ B(F) $ maps to the class of~$ \a $
in $ B(F(t)) $. This is unsatisfactory because $ b $
may have to depend on $ \a $, so this does not give a uniform
map that induces the inverse of the isomorphism in Proposition~\ref{susprop},
and it may fail completely if~$ F $ is finite.

If $ |F| \ge 4 $, then applying $ \sp b,c $ as in Section~\ref{results}
for $ b $ in $ F \cup \{\infty\} $ and $ c $ in $ F^\flat $
to~$ \Z[F(t)^\flat] $ maps~$ \Rft F(t) $ to $ \Rft F $. So we
make the following definition, which is similar to but different
from one in~\cite[\S3.3]{BdJGRY}.

\begin{definition}
For any field $ F $, we let $ \pptilde(F) = \Z[F^\flat] / \Rft F $.
\end{definition}

Then $ \twt F^* $ as defined in Section~\ref{results} is a quotient of $ \suswedge F^* $,
$ \pptilde(F) $ is a quotient of~$ \pp(F) $, and 
the map~$ \delta_F $ induces $ \partial_F : \Z[F^\flat] \to \twt F^* $
as introduced in Section~\ref{results}, mapping $ [x] $ to $ x \tw (1-x) $.
As~$ \partial_F(R_2(F)) = 0 $, we can define the modified
Bloch group $ \Btilde(F) = \ker(\partial_F) / \Rft F $
inside~$ \pptilde(F) $.
We then have the following analogue of Proposition~\ref{susprop}.

\begin{proposition} \label{barsusprop}
If $ |F| \ge 4 $, and $ t $ is a variable, then the map $ \Btilde(F) \to \Btilde(F(t)) $ is an isomorphism.
\end{proposition}

\begin{proof}
The class $ C_F $ in $ \pp(F) $ of any $ C_{F,c} $ is independent of $ c $ in $ F^\flat $
because $ |F| \ge 4 $.
One shows as in the proof of \cite[Theorem~3.23]{BdJGRY}
that $ \Btilde(F) = B(F)/\langle 3 C_F \rangle $
by replacing $ \pp(F)/\langle C_F \rangle $ with~$ \pp(F)/\langle 3 C_F \rangle $
and~$ \ol \pp (F) $ with~$ \pptilde (F) $ in the diagram preceding
the theorem of loc.~cit., using that $ 3 C_{F,c} $ is in $ \Rft F $.
(One can relax the assumption that $ F $ is infinite in loc.~cit.~to~$ |F| \ge 4 $ by following \cite[VI.5]{WeiKbook} instead of~\cite{susXkof}
for the set-up and the proof.)
The proposition then follows from Proposition~\ref{susprop} as
$ C_F $ in~$ B(F) $ maps to~$ C_{F(t)} $ in~$ B(F(t)) $.
\end{proof}

If $ |F| \ge 4 $, and $ \a $ in $ \Z[F(t)^\flat] $ satisfies~$ \partial_{F(t)}(\a) = 0 $,
then by Proposition~\ref{barsusprop} there exists~$ \bb $ in $ \Z[F^\flat] $
with~$ \partial_F(\bb ) = 0 $ and~$ \a - \bb $ in $ \Rft F(t) $.
Applying any $ \sp b,c $ we find~$ \sp b,c (\a) - \bb $ is in
$ \Rft F $, so that~$ \sp b,c (\a) $ is
in the kernel of $ \partial_F $.
(This also follows in an elementary way from~\eqref{spCD} below.)
Hence~$ \sp b,c $ always induces the inverse of the isomorphism
in Proposition~\ref{barsusprop}.

Now let $ A $ be any torsion-free additive group.
We noted in Section~\ref{results} that tensoring with $ A $ over
$ \Z $ is exact, hence
tensoring $ \partial_F : \Z[F^\flat] \to \twt F^* $ with $ A $
over $ \Z $ shows that $ \ker(\id_A \otimes \partial_F) = A \otimes_\Z \ker(\partial_F) $.
So if~$ |F| \ge 4 $, then from the isomorphism
\begin{equation*}
 \frac{\ker( \Z[F^\flat] \to \twt F^*)}{\Rft F }
\to
 \frac{\ker(\Z[F(t)^\flat] \to \twt F(t)^*)}{\Rft F(t) }
\end{equation*}
in Proposition~\ref{barsusprop} we obtain an isomorphism
\begin{equation} \label{explicit}
 \frac{\ker( A[F^\flat] \to A \otimes_\Z \twt F^*) }{A \Rft F }
\to
 \frac{\ker( A[F(t)^\flat] \to A \otimes_\Z \twt F(t)^*)}{A \Rft F(t) }
\,,
\end{equation}
for which the inverse is still given by (the map induced by) $ \sp b,c $ for any~$ b $ in~$ F \cup \{\infty\} $ and $ c $ in $ F^\flat $.
In order to simplify notation, we suppress $ \id_A $ from the
notation, writing simply $ \partial_F $ or $ \partial_{F(t)} $,
etc.

We now relate $ \sp b,c $, $ \partial_F $ and $ \partial_{F(t)} $.
From a local parameter $ \pi_b $ at~$ b $ of $ \P_F^1 $
we obtain a homomorphism from~$ F(t)^* $ to $ F^* $ by mapping
$ f(t) $ to $ \pi_b^{-\ord_b(f(t))} f(t) _{|t=b} $. It induces
a homomorphism~$ \spp b : \twt F(t)^* \to \twt F^* $ that depends
on the choice of~$ \pi_b $, but the diagram
\begin{equation} \label{spCD}
\begin{split}
\xymatrix{
 \Z[F(t)^\flat] \ar[r]^{\partial_{F(t)}} \ar[d]_{\sp b,c } & \twt F(t)^* \ar[d]^{\spp b }
\\
 \Z[F^\flat] \ar[r]^{\partial_F} & \twt F^*
}
\end{split}
\end{equation}
always commutes. (Note that $ \partial_F ( [c]+[1-c]) = 0 $.)

\begin{lemma} \label{endlemma}
Let $ \a $ be in $ A[F(t_1,\dots,t_n)^\flat] $ with
$ n \ge 1 $ and $ |F| \ge 4 $. Suppose that $ \partial(\a) $ is in $ A \otimes_\Z \twt F^* $.
Then~$ \a - \a' $ is in~$ A \Rft F(t_1,\dots,t_n) $ for any iterated specialisation $ \a' $ of $ \a $.
Moreover, $ \a' - \a'' $ is in~$ A \Rft F(t_1,\dots,t_m) $
if $ \a' $ and $ \a'' $ are iterated specialisations of
$ \a $ to $ F(t_1,\dots,t_m) $ with~$ 0 \le m < n $.
\end{lemma}

\begin{proof}
First assume $ n = 1 $.
If $ \a $ in $ A[F(t)^\flat] $ has $ \partial(\a) $
in $ A \otimes_\Z \twt F^* $, then $ \a - \sp b,c (\a) $ is in $ \ker(\partial) $
for any $ b $ in $ F \cup \{\infty\} $ and $ c $ in $ F^\flat $
by the commutativity of~\eqref{spCD}.
It follows from~\eqref{explicit} that $ \a - \sp b,c (\a) $ is in~$ A \Rft F(t) $
because $ \sp b,c $ induces the inverse isomorphism.
For $ n \ge 2 $ one applies this at every step of the iterated
specialisation and takes a telescoping sum of such differences $ \a - \sp b,c \a $.
For the last statement, one now has~$ \a' - \a'' $ in
$ A \Rft F(t_1,\dots,t_n) $, and iterated
specialisation to $ F(t_1,\dots,t_m) $ maps this to
$ A \Rft F(t_1,\dots,t_m) $ without changing $ \a' $ or $ \a'' $.
\end{proof}

\begin{example}
If $ \partial(\a) = 0 $ in Lemma~\ref{endlemma}, and $ P $ in
$ F^n $ is such that all functions involved in $ \a $ are defined at~$ P $ with value
not equal to 0 or 1, then 
$ \a - \a(P) $ is in~$ A \Rft F(t_1,\dots,t_n) $
and $ \partial(\a(P)) = 0 $.
So $ \a(P) $ represents the class in the left-hand side of~\eqref{explicit}
that maps the class of $ \a $ under iterations of the map in~\eqref{explicit}.
\end{example}

We need more ingredients for our proofs.
Let $ F $ be any field.
For each valuation~$ v : F^* \to \Z $ with residue field~$ \kappa_v $,
we define the homomorphism~$ T_v : \twt F^* \to \kappa_v^* $
by mapping
$ f \tw g $ to the class of~$ (-1)^{v(f) v(g)} f^{v(g)} g^{-v(f)} $.
It is well-known (and easily checked) that
it is trivial on the image of $ \partial $.
For two such valuations $ w $ and $ \ww $ we define 
$ T_{w,\ww} : \twt F^* \to \Z $ by mapping $ f \tw g $
to~$ w(f) \ww(g) - \ww(f) w(g) $.

If $ t_1,\dots,t_n $ are variables, and $ A $
is any additive group, we already noticed in Section~\ref{results} that $ A \otimes_\Z \twt F^* $ injects
into $ A \otimes_\Z \twt F(t_1,\dots,t_n)^* $.
If $ A $ is also torsion-free,
then inside~$ A \otimes_\Z \twt F(t_1,\dots,t_n) $ we have
\begin{equation*}
 A \otimes_\Z \twt F^* = \mybigcap _v \ker(\id_A \otimes T_v) \cap \mybigcap_{w, \ww} \ker(\id_A \otimes T_{w,\ww}) 
\end{equation*}
where $ v $ and the $ w \ne \ww $ run through the valuations associated to
the irreducibles of~$ F[t_1,\dots,t_n] $.
More precisely,
if we choose representatives of those irreducibles up to association,
then any element of  $ A \otimes_\Z \twt F(t_1,\dots,t_n)^* $
can be written as
\begin{equation} \label{betas}
  \b_1 + \b_2 + \b_3 
\end{equation}
with $ \b_1 $ a sum of terms $ a \otimes (p \tw q) $ with $ p,q $
among the chosen irreducibles, and~$ p \ne q $ as $ p \tw p = p \tw (-1) $,~$ \b_2 $ a sum of terms $ a' \otimes (p \tw c) $
with $ p $ a chosen irreducible and $ c $ in $ F^* $, and $ \b_3 $ in $  A \otimes_\Z \twt F^* $.
Then $ \id_A \otimes T_{w,\ww} $ for $ w $ and $ \ww $ the
valuations corresponding to~$ p \ne q $ maps $ a \otimes (p \tw q) + b \otimes( q \tw p) $
to $ a - b $ and is trivial on all other terms in~\eqref{betas}, which shows that~$ \b_1 $ is unique.
For the valuation $ v $ corresponding to $ p $,
note that~$ A \otimes_\Z F^* $ injects into~$ A \otimes_\Z \kappa_v^* $
because $ A $ is torsion-free.
As the contribution $ \sum_i a_{p,i} \otimes ( p \tw c_i) $ to
$ \b_2 $ for $ p $ is mapped under $ \id_A \otimes T_v $ to~$ \sum_i a_{p,i} \otimes c_i $,
and this map is trivial on all other terms in $ \b_2 $ as well
as on $ \b_3 $, we see $ \b_2 $ is unique as well. The same then holds for
$ \b_3 $.

Because each $ T_v $ is trivial on the image of $ \partial $,
we also have
\begin{equation} \label{imker}
\im(\id_A \otimes \partial_{F(t_1,\dots,t_n)} ) \cap \mybigcap _{w, \ww} \ker(\id_A \otimes T_{w,\ww})
\subseteq  A \otimes_\Z \twt F^*
\,.
\end{equation}

\begin{remark}
In the above we used $ p \tw p = p \tw (-1) $,
whereas in Suslin's construction one has~$ 2 (\sustimes p p ) = 0 $.
This detail, which corresponds to using $ R_{5,2} $ instead
of $ R_5 $, allows us to avoid a multiplication by 2 in the statements
of the propositions in Section~\ref{results}.
\end{remark}

For dealing with more than one variable in our proofs, we shall
use the following.

\begin{lemma} \label{buddinglemma}
Let $ A $ be an additive torsion-free group,
$ \tF $ an infinite subfield of the field $ F $, and $ t_1,\dots,t_n $ variables with~$ n \ge 2 $.
Suppose that $ \a $ in $ A[F(t_1,\dots,t_n)^\flat] $
under $ \partial $ is not mapped into~$  A \otimes_\Z \twt F^* $.
Then there exists $ P = (b_1,\dots, b_n) $ in $ \tF^n $, and an
index~$ l $ in~$ \{ 1,\dots,n \} $, such that
all functions involved in $ \a $ are defined at~$ P $ with value
not equal to 0 or 1, 
and $ \a(b_1,\dots,b_{l-1}, t_l, b_{l+1}, \dots, b_n) $
under $ \partial $ is not mapped into~$ A \otimes_\Z \twt F^* $.
\end{lemma}

\begin{proof}
If $ \partial(\a) $ is
in the kernel of $ \id_A \otimes T_{w,\ww} $ for all $ w \ne \ww $
associated to irreducibles in $ k[t_1,\dots,t_n] $
with $ \deg_{t_n} $ positive, then it is in $ A \otimes_\Z \twt F(t_1,\dots,t_{n-1})^* $
by~\eqref{imker}
with $ F $ and~$ F(t_1,\dots,t_n) $ replaced with~$ F(t_1,\dots,t_{n-1}) $ and $ F(t_1,\dots,t_{n-1})(t_n) $ respectively.
Therefore
replacing $ \a $ with $ \sp 0,c (\a) $, where we specialise $ t_n $ to 0 and~$ c $ is
in~$ F^\flat $, does not change the image under $ \partial $
but reduces the number of variables in $ \a $.
Continuing this way, we find an 
index $ l $ for which there are two non-associate irreducible
$ p $ and $ q $ in~$ F[t_1,\dots,t_l] $ with $ \deg{t_l} $ positive and for which the associated valuations
$ w $ and $ \ww $ give $ (\id_A \otimes T_{w,\ww}) (\partial(\a)) \ne 0 $.
We emphasise this also applies to the original~$ \partial(\a) $ as it
did not change under the specialisations.

Returning to our original $ \a $, we shall assume $ l = n $ for
notational simplicity. 
Choose finitely many non-associate irreducibles $ p_j $ in $ F[t_1,\dots,t_n] $
that represent all the irreducible factors of the numerators
and denominators of $ f $ and $ 1-f $ for the elements $ [f] $
that are used in $ \a $, so that we can express $ \partial(\a) $
as in~\eqref{betas} using only those $ p_j $ and elements of~$ F^* $.

Using that $ \tF $ is infinite, now choose $ P = (b_1,\dots,b_n) $ in $ \tF^n $ such that:
\begin{itemize}
\item
$ p_j(P) \ne 0 $ for all $ p_j $;

\item
$ \deg_{t_n}(p_j(b_1,\dots, b_{n-1},t_n)) = \deg_{t_n}(p_j) $ if $ \deg_{t_n}(p_j) $ is positive;

\item
for each $ p_j \ne p_{j'} $ with $ \deg_{t_n}(p_j) $ and $ \deg_{t_n}(p_{j'}) $
both positive, the resultant
$ \Res_{t_n}(p_j, p_{j'}) $, which is non-zero in $ F[t_1,\dots,t_{n-1}] $,
does not evaluate to zero at $ b_1,\dots,b_{n-1} $.
\end{itemize}
The first condition implies that all functions involved in $ \a $
are defined at $ P $ with value not equal to 0 or 1.
The second condition means that computing the resultant
in the third commutes with specialising~$ t_1,\dots, t_{n-1} $.
Combined, they imply that for $ p_j \ne p_{j'} $ with $ \deg(t_n)(p_j) $
and $ \deg_{t_n}(p_{j'}) $ both positive,
$ p_j(b_1,\dots,b_{n-1}, t_n) $ and~$ p_{j'}(b_1,\dots,b_{n-1}, t_n) $
are not units and relatively prime in $ F[t_n] $.
So, for every~$ p_j $ with $ \deg_{t_n}(p_j) > 0 $
there is an irreducible $ \ol p_j $ in $ F[t_n] $ dividing~$ p_j(b_1,\dots,b_{n-1}, t_n) $ but
no $ p_{j'}(b_1,\dots,b_{n-1},t_n) $ if $ p_{j'} \ne p_j $.

Let $ \ol p $ and $ \ol q $ in $ F[t_n] $ be such irreducibles for
the irreducibles $ p $ and $ q $ in~$ F[t_1,\dots,t_n] $ that correspond to
the~$ w $ and $ \ww $ we found above.
Let $ m_p $ be the multiplicity of $ \ol p $ in~$ p(b_1,\dots,b_{n-1},t_n) $
and $ m_q $ the multiplicity of $ \ol q $ in~$ q(b_1,\dots,b_{n-1},t_n) $.
If $ w' $ and~$ \ww' $ are the valuations of
$ F(t_n) $ corresponding to $ \ol p $ and $ \ol q $, then
\begin{equation*}
 (\id_A \otimes T_{w',\ww'}) (\partial(\a(b_1,\dots,b_{n-1},t_n))) = m_p m_q (\id_A \otimes T_{w,\ww}) (\partial(\a)) \ne 0 
\end{equation*}
in $ A $ because it is torsion-free.
Hence $ \partial(\a(b_1,\dots,b_{n-1},t_n)) $
is not in $ A \otimes_\Z \twt F^* $.
\end{proof}

We can now give the proofs of the propositions in Section~\ref{results}.
Except for Proposition~\ref{conjpropC}, the proof proceeds
by showing that the dilogarithm involved being constant on $ \a $
implies~$ (\id_A \otimes \partial) (\a) $ is in $ A \otimes_\Z \twt k^* $
in the case of one variable by using~\eqref{imker}, then invoking Lemma~\ref{buddinglemma}
to deduce this for two or more variables, and finally applying Lemma~\ref{endlemma}.
The first step depends on the specific dilogarithm involved, the last two
steps barely or not at all.
Proposition~\ref{conjpropC} is proved in a different way
because one cannot specialise a complex variable and its conjugate
independently.

For Proposition~\ref{funprop} 
we note $ \wirt z D(z) =  \frac1{2i} ( \log|z| \wirt z \log(1-z) - \log|1-z| \wirt z \log(z) ) $.
So, for $ \a $ in $ A[k(t)^\flat] $, we can obtain $ \wirt t D(\a) $ by
first computing $ \partial(\a) $ in $ A \otimes_\Z \twt k(z)^* $,
and then applying to it the map that maps $ a \otimes (f \tw g) $
to
\begin{equation*}
\frac a{2i} (\log|f| \wirt t \log(g) - \log|g| \wirt t \log(f)) = \frac a{2i} (\log|f| \frac{g'}{g} - \log|g| \frac{f'}{f})
\,.
\end{equation*}
In order to get cleaner expressions, we shall use the scaled version $ 2 i \wirt t D(\a) $ instead.

\begin{proof}[Proof of Proposition~\ref{funprop}.]
We first let $ n = 1 $ and write $ t $ for $ t_1 $.
Choosing the monic irreducibles in $ k[t] $, we write~$ \partial(\a) = \b_1 + \b_2 + \b_3 $ as in~\eqref{betas}.
Let~$ z = t - \c $ with $ \c $ a root in~$ \C $ of (exactly)
one monic irreducible~$ p(t) $.
The contribution to~$ 2 i \wirt t D(\a) $ of a term
in $ \partial(\a) $ has its expansion around~$ \c $ in $ \C[[z,\ol z]] $,
except for~$ \sum_j a_j \otimes (p(t) \tw q_j(t)) $
with $ q_j(t) \ne p(t) $ monic irreducible, from $ \b_1 $, and $ \sum_j a_j' \otimes (p(t) \tw c_j) $,
from~$ \b_2 $. They give
\begin{equation*}
\sum_j a_j \log|p(t)|\frac{q_j'(t)}{q_j(t)}
-
\sum_j a_j \log|q_j(t)| \frac{p'(t)}{p(t)}
-
\sum_j a_j' \log|c_j| \frac{p'(t)}{p(t)}
\,.
\end{equation*}
The expansion of the first term here is the sum of
$ \sum_j a_j \log|z| \frac{q_j'(z+\c)}{q_j(z+\c)} $ and an element
of~$ \C[[z, \ol z]] $,
the other two terms have expansions in $ z^{-1} \C[[z, \ol z]] $.
Because~$ \log|z| \, \C[[z]] \cap z^{-1} \C[[z,\ol z]] = 0 $,
and~$ \wirt t D(\a) = 0 $ we see that~$ \sum_j a_j \frac{q_j'(z+\c)}{q_j(z+\c)} = 0 $,
hence all $ a_j $ are zero because the $ q_j(t) $ have no common
zeroes.
Applying this to all $ p(t) $ gives~$ \b_1 = 0 $, so that
$ \partial(\a) $ is in $ A \otimes_\Z \twt k^* $ by~\eqref{imker}.

For $ n > 1 $ we proceed by contradiction.
If $ \partial(\a) $ is
not in $ A \otimes_\Z \twt k^* $, then by Lemma~\ref{buddinglemma} there is
an index~$ l $ and a point $ P = (b_1,\dots,b_n) $ in $ k^n $
such that on one hand specialising all $ t_j $ for $ j \ne l $
to $ b_j $ in $ \a $ gives an element~$ \a' $ of~$ A[k(t_l)^\flat] $
with~$ \partial(\a') $ not in $ A \otimes_\Z \twt k^* $, but on the
other hand 
$ D(\a') $ is constant and defined on
a Zariski open part
of~$ \{b_1\} \times \dots \times \{b_{l-1}\} \times k \times \{b_{l+1}\} \times \dots \times \{b_n\} $
by our choice of $ P $ as it is the restriction of $ D(\a) $
to this subset of $ k^n $.
This contradicts the case~$ n = 1 $.

Because $ \partial(\a) $ is in $ A \otimes_\Z \twt k^* $ for
all $ n \ge 1 $ by the
above, we can now apply Lemma~\ref{endlemma} for the statement
about iterated specialisations.
As explained before Example~\ref{introexample},
using a point at which all functions in $ \a $ are
defined with value not equal to 0 or 1 is just a special case
of this.
\end{proof}

\begin{remark} \label{sharpDrem}
Let $ A = \Z $, and take $ f $ in $ k(t_1,\dots,t_n) \setminus k $ such
that $ c f $ for $ c $ in~$ k^* $ is never a square. If~$ \a = [f] + [f^{-1}] $,
then $ \a - \a' $ is in $ R_2 $ for any $ n $-step iterated specialisation
$ \a' $, but it is not in $ R_5 $ as it is not in the kernel of Suslin's map~$ \delta_{k(t_1,\dots,t_n)} $ 
(see the proof of Proposition~\ref{barsusprop}).
Thus Proposition~\ref{funprop} seems sharp.
Similar considerations apply to Propositions~\ref{Rfunprop} and
Proposition~\ref{Coleprop}.
\end{remark}

\begin{proof}[Proof of Proposition~\ref{conjpropR}.]
We first take $ n = 1 $ and again write $ t $ for $ t_1 $.
We compute a scaled version of~$ \frac{\dd}{\dd t} D(\a) $
using the observation before the proof of Proposition~\ref{funprop}. Any term
$ a \otimes (f(t) \tw g(t)) $ in~$ \partial(\a) $ then contributes
\begin{equation} \label{realcont}
  a \log(f(t) \ol f(t)) \Bigl(\frac{g'(t)}{g(t)} - \frac{\ol g'(t)}{\ol g(t)} \Bigr)
- a \log(g(t) \ol g(t)) \Bigl(\frac{f'(t)}{f(t)} - \frac{\ol f'(t)}{\ol f(t)} \Bigr)
\end{equation}
because $ t $ is real.
Therefore $ \frac{\dd}{\dd t} D(\a) $
is the restriction to a Zariski open part of $ \R $ of a 
holomorphic function on the complement in $ \C $ of several vertically
downward cuts that make all logarithms involved single-valued and remove
any poles.
Because it is identically zero on an open interval in $ \R $, the
same holds on this complement in $ \C $. Our argument below will
be on the monodromy of~\eqref{realcont} around roots of $ f(t) \ol f(t) $ and~$ g(t) \ol g(t) $.

We choose the monic irreducibles in $ k[t] $ and write
$ \partial(\a) = \b_1 + \b_2 + \b_3 $ as in~\eqref{betas}.
If $ p(t) $ and~$ \ol p(t) $ are associate monic irreducibles then they are equal, so
this decomposition is compatible with complex conjugation.
For two monic irreducibles~$ p(t) \ne q(t) $ we consider,
because of~\eqref{realcont}, the terms
of $ \b_1 $ that contain exactly two of $ p(t), q(t), \ol p(t) $ and $\ol q(t) $,
which can give the following types.

\mybullet
$ a \otimes ( p(t) \tw \ol p(t) ) $ with $ \ol p(t) \ne p(t) $.
It contributes $ 2 a \log(p(t) \ol p(t)) \bigl(\frac{\ol p'(t)}{\ol p(t)} - \frac{p'(t)}{p(t)} \bigr) $
in~\eqref{realcont}, with monodromy around a root of $ p(t) $ proportional
to $ 2 a \bigl(\frac{\ol p'(t)}{\ol p(t)} - \frac{p'(t)}{p(t)} \bigr) $.
As this cannot can from any other term in $ \partial(\a) $, it
follows that $ a = 0 $.

\mybullet
$ a \otimes (p(t) \tw q(t)) $
with $ \ol p(t) = p(t) $ and $ \ol q(t) = q(t) $, which is invariant
under complex conjugation.

\mybullet
$ a \otimes (p(t) \tw q(t)) + a' \otimes (p(t) \tw \ol q(t)) $
where $ \ol p(t) = p(t) $ but $ \ol q(t) \ne q(t) $. From the
monodromy around a root of $ p(t) $ we find $ a = a' $, so that
it is invariant under complex conjugation.

\mybullet
$ a \otimes (p(t) \tw q(t)) + a' \otimes (p(t) \tw \ol q(t)) + b \otimes (\ol p(t) \tw q(t)) + b' \otimes (\ol p(t) \tw \ol q(t)) $
with~$ \ol p(t) \ne p(t) $ and $ \ol q(t) \ne q(t) $.
From the monodromy around a root of $ p(t) $ we find that
$ a - a' + b - b' = 0 $ and from the monodromy around a root
of $ q(t) $ that~$ a + a' - b - b' = 0 $. Therefore $ a = b' $ and $ a' = b $,
and the sum is invariant.

This shows that $ \ol \b_1 = \b_1 $.
Then $ \b_2 - \ol \b_2 + \b_3 - \ol \b_3 $ is the decomposition
of $  \partial(\a - \ol \a) $, and it follows from~\eqref{imker}
and the discussion of~\eqref{betas}
that $ \ol \b_2 = \b_2 $.

Now take $ n > 1 $.
When choosing irreducibles up to association in $ k[t_1,\dots,t_n] $,
we can again assume that if~$ p $ and $ \ol p $ are associate then $ p = \ol p $
by scaling a coefficient to~1.
Then the decomposition in~\eqref{betas} is again compatible with
complex conjugation, so if $ \partial(\a) = \b_1 + \b_2 + \b_3 $ then
$ \partial(\a - \ol \a) = (\b_1 - \ol \b_1) + (\b_2 - \ol \b_2) + (\b_3 - \ol \b_3) $.
If $ \b_1 = \ol \b_1 $ then $ \b_2 = \ol \b_2 $ by~\eqref{imker} and we are done.

To rule out $ \b_1 \ne \ol \b_1 $ we proceed by contradiction.
If $ \b_1 - \ol \b_1 \ne 0 $ then
$ \partial(\a - \ol \a) $ is not in $ A \otimes_\Z \twt k^* $.
Using Lemma~\ref{buddinglemma} we can find a point $ P = (b_1,\dots,b_n) $ in~$ k^n \cap \R^n $ and
an index $ l $ such that, with $ \a_1 $ the partial specialisation of
$ \a $ obtained by specialising $ t_j $ to $ b_j $ for $ j \ne l $ and
$ \a_1' $ similarly obtained from~$ \ol \a $,~$ \partial(\a_1 - \a_1') $ is not in $ A \otimes_\Z \twt k^*  $ but $ D(\a_1-\a_1') $
is constant as it is the restriction of~$ D(\a - \ol \a) = 2 D(\a) $.
From~$ P = \ol P $ we have~$ \a_1' = \ol \a_1 $, so that $ D(\a_1) = \frac 12 D(\a_1 - \ol \a_1)  $ is constant
on the Zariski open part of $ \R $ where it is defined but $ \partial(\a_1 - \ol \a_1) $
is not in $ A \otimes_\Z \twt k^* $, which contradicts the case
$ n = 1 $.

We now know that $ \partial(\a - \ol \a) = \b_3 - \ol \b_3 $ is in $ A \otimes_\Z \twt k^* $
for all $ n \ge 1 $, so, as in the proof of Proposition~\ref{funprop},
the statement about iterated specialisations follows from Lemma~\ref{endlemma}.
\end{proof}

We move on to the proof of Proposition~\ref{conjpropC}.
For $ k $ a subfield of $ \C $, the ring of functions on $ \C^n $ given by
elements of~$ k[z_1,\ol z_1, \dots, z_n, \ol z_n] $ 
is isomorphic to the polynomial ring $ k[s_1,\dots, s_{2n}] $
by mapping $ s_{2j-1} $ to $ z_j $ and $ s_{2j} $ to~$ \ol z_j $.
(This is an injection as it is injective on $ k[z_1,\dots,z_n] $,
and if any $ \ol z_j $ were to occur in an element of the kernel then one gets
a contradiction by taking one with $ \deg_{\ol z_j} $ minimal
and applying $ \wirt \ol z_j $ to~it.)

\begin{proof}[Proof of Proposition~\ref{conjpropC}.]
If $ k $ contains $ i $, then letting $ z_j $ correspond to $ s_{2j-1} + i s_{2j} $
and~$ \ol z_j $ to $ s_{2j-1} - i s_{2j} $ gives
an isomorphism $ k[s_1,\dots,s_{2n}] \to k[z_1,\ol z_1,\dots,z_n, \ol z_n] $
that is compatible with the complex conjugations on either side.
From Proposition~\ref{conjpropR} we know that~$ \partial(\a) $ is the sum of
an element in~$ A \otimes_\Z \twt k(z_1, \ol z_1,\dots,z_n, \ol z_n)^* $ that is invariant under complex conjugation,
and an element in~$ A \otimes_\Z \twt k^* $.
So~$ \partial( \a - \ol \a ) $ is in~$ A \otimes_\Z \twt k^* $.

If $ i $ is not in $ k $, then
for the isomorphism $ k(z_1,\ol z_1,\dots,z_n, \ol z_n) \to k(s_1,\dots,s_{2n}) $
in the discussion just before this proof,
the image of $ \partial( \a - \ol \a ) $ in $ \twt k(s_1,\dots,s_{2n})^* $ equals $ \b_1 + \b_2 + \b_3 $
as in~\eqref{betas} for some choice of irreducibles of $ k[s_1,\dots,s_{2n}] $ up to
association.
The isomorphism is compatible with extending $ k $ to $ k(i) $,
and one can get a corresponding decomposition $ \b_1' + \b_2' + \b_3' $
from $ \b_1 + \b_2 + \b_3 $
by making a choice of irreducibles up to association for $ k(i)[s_1,\dots,s_{2n}] $,
factorising the chosen irreducibles for~$ k[s_1,\dots,s_{2n}] $
into those, and some rewriting.
We already know that~$ \b_1' = \b_2' = 0 $ and considering
how an irreducible in $ k[s_1,\dots,s_{2n}] $ can factorise after extending 
$ k $ to $ k(i) $ one sees that $ \b_1 = 0 $. It then follows
from~\eqref{imker} that~$ \b_2 = 0 $ as well, which proves what
we want.

For the point $ P = (b_1,\dots,b_n) $ in $ k^n $, we notice that the assumptions are such that for
every $ [f] $ involved in~$ \a $, the irreducibles in the numerators and
denominator of $ f $ and $ 1-f $ in lowest terms are non-zero
at~$ P $, hence $ \a(P) - \ol{\a(P)} $ is obtained by evaluating the functions
at $ P $.
Under the isomorphism of the field with $ k(s_1,\dots,s_{2n}) $
in the discussion just before the beginning of this proof, this corresponds to an iterated
specialisation that specialises all $ s_{2j-1} $ to $ b_j $
and all~$ s_{2j} $ to~$ \ol b_j $.
So it follows from Lemma~\ref{endlemma} and this isomorphism
that~$ \a - \ol \a - \a(P) + \ol{\a(P)} $ is in $ A \Rft k(z_1,\ol z_1,\dots,z_n,\ol z_n) $.
\end{proof}

\begin{remark} \label{Dbarsharprem}
In the context of Proposition~\ref{conjpropR},
let
\begin{equation*}
 R_{cc}(k(t_1,\dots,t_n)) = \langle [f] + [\ol f] \text{ with } f \text{ in } k(t_1,\dots,t_n)^\flat \rangle 
\end{equation*}
in $ \Z[k(t_1,\dots,t_n)^\flat] $ be
the subgroup of relations corresponding to the functional equation~$ D(z) + D( \ol z) = 0 $
involving complex conjugation.
From the proposition we see that if~$ D(\a) $ is constant, then~$ 2 (\a-\a') $ is in $ A R_{5,2} + A R_{cc} $ for
any iterated specialisation~$ \a' $ of $ \a $ that commutes with
complex conjugation (see Remark~\ref{conjremR}).
But $ \a-\a' $ may not be in there. For example, let $ n = 1 $,  $ A = \Z $, $ k $ any subfield
of~$ \R $, and~$ \a = [t] $. Then $ D(\a) = 0 $ on $ \R^\flat $,
$ \a' = [c] $, 0 or~$ [c] + [1-c] $ or~$ -[c] - [1-c] $ for
some $ c $ in $ k^\flat $, 
so $ \a - \a' $ is not in $ R_{5,2} + R_{cc} $
as $ \partial(R_{5,2} + R_{cc}) = \partial(R_{cc}) $ is mapped
to $ 2 \Z $ under any $ T_{w,\ww} $ for valuations $ w $ and~$ \ww $ corresponding to irreducibles of $ k[t] $
since~$ \ol f = f $, but $ \partial(\a - \a') $  is mapped to
1 if we let $ w $ correspond to $ t $ and $ \ww $ to~$ t-1 $.
Thus Proposition~\ref{conjpropR} and this consequence involving
$ R_{cc} $ both seem sharp.

Similar considerations apply to Proposition~\ref{conjpropC} with
the same $ n $, $ A $ and $ k $, and~$ \a = [z+\ol z] $.
Then~$ \a(P) = [c] $ for $ c $ in $ k^\flat $, and $ \a - \a(P) $ is not in
$ R_{5,2} + R_{cc} $ as one can see using the isomorphism~$ k(z,\ol z) \simeq k(s_1,s_2) $ and considering valuations corresponding
to irreducibles in $ k[s_1,s_2] $ that are symmetric in $ s_1 $
and~$ s_2 $.
\end{remark}

Before proving Proposition~\ref{Rfunprop}, we discuss
some functional equations of~$ \L $ and~$ \RL $.
As stated on~\cite[p.23]{Zag07}, for $ \L $ we have
$ \L(x) + \L(1-x) = \l $ and $ \L(x) + \L(x^{-1}) = - \l $  for $ x $ in $ \P_\R^1  $,
and
\begin{equation*}
 \L\Bigl(\frac{1-x}{1-xy}\Bigr) + \L(x) + \L(1-xy) + \L(y) + \L\Bigl(\frac{1-y}{1-xy}\Bigr) = 0  
\end{equation*}
for $ x $ and $ y $ in $ \P_\R^1 \times \P_\R^1 \setminus \{ (1,1) , (\infty,0) , (0,\infty) \} $,
so that it does not vanish on~$ R_2 $.
(The arguments $ f_1,\dots,f_5 $ in this order satisfy $ f_{i-1} f_{i+1} + f_i = 1 $
with indices modulo~5.)
But the map induced by $ \RL(x) = \L(x) - \l $ vanishes even
on the larger subgroup of~$ \Z[\R \cup \{\infty\}] $ in~\eqref{specsubgroup} because we have (cf.~\cite{susXkof})
\begin{itemize}
\item
$ \RL(x) + \RL(1-x) + \l = 0 $ for $ x $ in $ \R^\flat $;

\item
$  \RL(x) + \RL(x^{-1}) = 0 $ for $ x $ in $ \R^\flat $;

\item
$  \RL(x)-\RL(y)+\RL(\frac{y}{x})+\RL(\frac{1-x}{1-y})-\RL(\frac{1-x^{-1}}{1-y^{-1}}) = 0 $
for $ x \ne y $ in $ \R^\flat $.
\end{itemize}
The last identity can be deduced from those for $ \L $ by using
\begin{alignat*}{1}
& \phantom{=}\,\,\,
[x^{-1}] - [y] + [xy] + \Bigl[ \frac{1-x^{-1}}{1-y} \Bigr] - \Bigl[ \frac{1-x}{1-y^{-1}} \Bigr]
\\
& =
- [x] - [y] - [1-xy] - \Bigl[ \frac{1-x}{1-xy} \Bigr] - \Bigl[ \frac{1-y}{1-xy} \Bigr]
+ ( [x] + [x^{-1}])
+ ( [xy] + [1-xy])
\\
& \phantom{=}\,\,\,
+ \Bigl( \Bigl[ \frac{1-x^{-1}}{1-y} \Bigr] +  \Bigl[ \frac{1-y}{1-x^{-1}} \Bigr]\Bigr)
- \Bigl( \Bigl[ \frac{1-y}{1-x^{-1}} \Bigr] +  \Bigl[ 1-\frac{1-y}{1-x^{-1}} \Bigr]\Bigr)
+ \Bigl( \Bigl[ \frac{1-xy}{1-x} \Bigr]     +  \Bigl[ \frac{1-x}{1-xy} \Bigr]\Bigr)
\\
& \phantom{=}\,\,\,
- \Bigl( \Bigl[ \frac{1-x}{1-y^{-1}} \Bigr] + \Bigl[ 1- \frac{1-x}{1-y^{-1}} \Bigr] \Bigr)
+ \Bigl( \Bigl[ \frac{1-xy}{1-y}  \Bigr]    + \Bigl[ \frac{1-y}{1-xy} \Bigr] \Bigr)
\end{alignat*}
because
$ 1-\frac{1-y}{1-x^{-1}} = \frac{1-xy}{1-x} $
and $ 1 - \frac{1-x}{1-y^{-1}} = \frac{1-xy}{1-y} $, and replacing
$ x $ with $ x^{-1} $.

\begin{proof}[Proof of Proposition~\ref{Rfunprop}.]
Locally on $ \R^\flat $ we can lift $ \RL(x) $ from $ \R/\frac{\pi^2}2 \Z $ to $ \R $ and then differentiate,
which gives~$ \RL'(x) = \frac12 \log|x| \, \dd \log|1-x| - \frac12 \log|1-x| \, \dd\log|x| $.
Applying the same principle to any $ \RL(\a) $ for~$ \a $ in $ A[k(t)^\flat] $,
with values in $ \C/\frac{\pi^2}2 A $, we get the result also by
first computing $ \partial(\a) $ in $ A \otimes_\Z \twt k(t)^* $,
and then applying to it the map that maps $ a \otimes (f \tw g) $
to
$\frac a4 (\log(f^2) \frac{g'}{g} - \log(g^2) \frac{f'}{f})$
because $ f(t) $ and~$ g(t) $ take values in~$ \R $.
For $ n = 1 $, the statement about $ \partial(\a) $ being in
$ A \otimes_\Z \twt k^* $ 
is then proved using a monodromy argument on the complement in $ \C $
of suitable cuts, similar to (but simpler than) what was given
in the proof of Proposition~\ref{conjpropR}.
The corresponding statement for $ n > 1 $,
and the final statement of the proposition, are
again proved using Lemma~\ref{buddinglemma} and Lemma~\ref{endlemma}.
\end{proof}

Before proving Proposition~\ref{Coleprop}, we discuss that
$ \Dp $ vanishes on $ R_{5,2} $.
From \cite[Proposition~6.4]{Col82} we have~$ \Dp(z) + \Dp(z^{-1}) = 0 $,
so that vanishes on $ R_2 $.
(It also satisfies $ \Dp(z) + \Dp(1-z) = 0 $ but we shall not
need this.)
Also, from \cite[Corollary~6.5b]{Col82} we have, with corrected signs,
that
\begin{equation*}
\Dp(xy) = \Dp(x) + \Dp(y) + \Dp\Bigl(\frac x {x-1} (1-y) \Bigr) + \Dp\Bigl(\frac y {y-1} (1-x) \Bigr)
\,.
\end{equation*}
Replacing $ x $ with $ x^{-1} $ and using the first of the functional
equations above gives that $ \Dp $ vanishes also on~$ R_5 $.

\begin{proof}[Proof of Proposition~\ref{Coleprop}.]
Let $ n = 1 $ and write $ t $ for $ t_1 $.
Because
\begin{equation*}
 \dd \Dp(z) = \frac12 \log_p(z) \, \dd \log_p(1-z) - \frac12 \log_p(1-z) \, \dd \log_p(z) 
\end{equation*}
we can compute $ \frac{\dd}{\dd z}\Dp(\a) $ by first computing $ \partial(\a) $,
and then mapping $ a \otimes (f \tw g) $ to~$ \frac a2 ( \log_p(f) \frac{g'}g - \log_p(g) \frac{f'}f ) $.
We write~$ \partial(\a) = \b_1 + \b_2 + \b_3 $ as in~\eqref{betas},
choosing the monic irreducibles in $ k[t] $.
Fix a root $ \c $ in~$ \C_p $ of such a monic~$ p(t) $, and let $ z = t - \c $.
The contribution to $ 2 \Dp'(\a) $ of a term
in $ \partial(\a) $ has its expansion around~$ \c $ in $ \C_p[[z]] $,
except for the terms~$ \sum_j a_j \otimes (p(t) \tw q_j(t)) $
with $ q_j(t) \ne p(t) $ monic irreducible, from~$ \b_1 $, and $ \sum_j a_j' \otimes (p(t) \tw c_j) $,
from~$ \b_2 $.
They give
\begin{equation*}
\sum_j a_j \log_p(p(t))\frac{q_j'(t)}{q_j(t)}
-
\sum_j a_j \log_p(q_j(t)) \frac{p'(t)}{p(t)}
-
\sum_j a_j' \log_p(c_j) \frac{p'(t)}{p(t)}
\,.
\end{equation*}
The expansion of the first term is the sum of
$ \sum_j a_j \frac{q_j'(t)}{q_j(t)} \log(z) $ and an element
of~$ \C[[z]] $,
the other two have expansions in $ z^{-1} \C[[z]] $.
As~$ \log_p(z) $ has derivative $ z^{-1} $, it is not
in $ \C_p((z)) $, so $ \log_p(z) \C_p[[z]] \cap z^{-1} \C_p[[z]] = 0 $
(cf.~\cite[p.182]{Col82}). It then
follows from $ 2 \Dp'(\a) = 0 $ that~$ \sum_j a_j \frac{q_j'(z+\c)}{q_j(z+\c)} = 0 $,
hence all $ a_j $ are zero because the~$ q_j(t) $ have no common
zeroes.
Applying this to all $ p(t) $ gives~$ \b_1 = 0 $, so that
$ \partial(\a) $ is in~$ A \otimes_\Z \twt k^* $ by~\eqref{imker}.

For $ n \ge 2 $ one again deduces that $ \partial(\a) $ is in $ A \otimes_\Z \twt k^* $ from
the case $ n = 1 $ by using Lemma~\ref{buddinglemma},
and Lemma~\ref{endlemma} gives the result on iterated specialisations.

For the dependency on the branch, first note that
$ \log_p^\circ(z) - \log_p(z) = \Delta v_p(z) $
for~$ \Delta = \log_p^\circ(p) - \log_p(p) $.
Therefore $ v_p(f) \log_p(g) - v_p(g) \log_p(f) $
for non-zero $ f $ and~$ g $ is
independent of the branch of the logarithm used in~it.

We also have $ \Lip^\circ(z) - \Lip(z) =  - \frac \Delta2 v_p(1-z) (\log_p^\circ(z) + \log_p(z)) $
by \cite[Proposition~2.6]{BdJ03}.
So for the branches of $ \Dp(z) = \Lip(z) + \frac12 \log_p(z) \log_p(1-z) $ we
have
\begin{alignat*}{1}
& \phantom{=}\,\,\,
\Dp^\circ(z) - \Dp(z)
\\
& =
 \Lip^\circ(z) -  \Lip(z) + \frac12 \left( \log_p^\circ(z) \log_p^\circ(1-z) - \log_p(z) \log_p(1-z) \right)
\\
& =
\frac\Delta2 \bigl(
 - v_p(1-z)  \left(\log_p^\circ(z) + \log_p(z)\right)
 + v_p(z) \log_p^\circ(1-z)
 + \log_p(z)v_p(1-z)
 \bigr)
\\
& =
\frac\Delta2 \left( v_p(z) \log_p(1-z) -  v_p(1-z) \log_p(z) \right)
\,.
\end{alignat*}
It follows that~$ \Dp^\circ(\a) - \Dp(\a) $, for any $ \a $ in $ A[k(t_1,\dots,t_n)^\flat] $,
can be directly computed from~$ \partial(\a) $ in $ A \otimes_\Z \twt k(t_1,\dots,t_n)^* $
by mapping $ a \otimes ( f \tw g) $ to
\begin{equation*}
a \frac{\log_p^\circ(p) - \log_p(p)}2 \left(v_p(f) \log_p(g) - v_p(g) \log_p(f)\right)
\,.
\end{equation*}
If $ \partial(\a) $ is in $ A \otimes_\Z \twt k^* $, this is
obviously a constant.
\end{proof}

\bibliographystyle{plain}
\bibliography{References}

\begin{thebibliography}{10}

\bibitem{BdJ03}
A.~Besser and R.~de~Jeu.
\newblock The syntomic regulator for the $ {K} $-theory of fields.
\newblock {\em Annales Scientifiques de l'\'Ecole Normale Sup\'erieure},
  36(6):867--924, 2003.

\bibitem{bl00}
S.~Bloch.
\newblock {\em Higher regulators, algebraic ${K}$-theory, and zeta functions of
  elliptic curves}.
\newblock Manuscript (`Irvine notes' 1978). Published as CRM Monograph Series,
  vol.~11 (American Mathematical Society, Providence, RI, 2000).

\bibitem{BdJGRY}
D.~Burns, R.~de~Jeu, H.~Gangl, A.~Rahm, and D.~Yasaki.
\newblock Hyperbolic tessellations and generators of $ {K}_3 $ for imaginary
  quadratic fields.
\newblock Submitted. Preprint available from http://arxiv.org/abs/1909.09091,
  2019.

\bibitem{Col82}
R.~Coleman.
\newblock Dilogarithms, regulators, and $p$-adic {$L$}-functions.
\newblock {\em Invent. Math.}, 69:171--208, 1982.

\bibitem{frenkel-szenes-1993}
E.~Frenkel and A.~Szenes.
\newblock Crystal bases, dilogarithm identities and torsion in algebraic
  ${K}$-groups.
\newblock {\em J. Amer. Math. Soc.}, 8(3):629--664, 1995.

\bibitem{gonXpam}
A.~B. Goncharov.
\newblock Polylogarithms and motivic {G}alois groups.
\newblock In {\em Motives (Seattle, WA, 1991)}, volume~55 of {\em Proc. Sympos.
  Pure Math.}, pages 43--96. Amer. Math. Soc., Providence, RI, 1994.

\bibitem{licht1989}
S.~Lichtenbaum.
\newblock Groups related to scissors-congruence groups.
\newblock In {\em Algebraic $K$-theory and algebraic number theory (Honolulu,
  HI, 1987)}, Contemp. Math., 83, pages 151--157. Amer. Math. Soc., Providence,
  RI, 1989.

\bibitem{parry-sah-1983}
W.~Parry and C.-H. Sah.
\newblock Third homology of ${SL}_2({\R})$ made discrete.
\newblock {\em J. Pure Appl. Algebra}, 30(2):181--209, 1983.

\bibitem{souderes18}
I.~Soud\`eres.
\newblock {\'E}quations fonctionnelles du dilogarithme.
\newblock {\em Annales de l'Institut Fourier}, 68(1):151--169, 2018.

\bibitem{sus:aka}
A.~A. Suslin.
\newblock Algebraic $ {K} $-theory of {F}ields.
\newblock In {\em Proceedings of the {I}nternational {C}ongress of
  {M}athematicians}, pages 222--243. 1986.

\bibitem{susXkof}
A.~A. Suslin.
\newblock {$K\sb 3$} of a field, and the {B}loch group (in {R}ussian).
\newblock {\em Trudy Mat. Inst. Steklov.}, 183:180--199, 229, 1990.
\newblock Translated in Proc.\ Steklov Inst.\ Math.\ {\bf 1991}, no.\ 4,
  217--239, Galois theory, rings, algebraic groups and their applications.

\bibitem{WeiKbook}
C.~Weibel.
\newblock {\em The {$K$}-book: an introduction to algebraic $K$-theory}, volume
  145 of {\em Graduate Studies in Mathematics}.
\newblock American Mathematical Society, Providence, RI, 2013.

\bibitem{Woj91}
Z.~Wojtkowiak.
\newblock A note on functional equations of the $p$-adic polylogarithms.
\newblock {\em Bull. Soc. Math. France}, 119(3):343--370, 1991.

\bibitem{Zag07}
D.~Zagier.
\newblock The dilogarithm function.
\newblock In {\em Frontiers in number theory, physics, and geometry. {II}},
  pages 3--65. Springer, Berlin, 2007.

\end{thebibliography}

\end{document}